\documentclass[11pt]{amsart}
\usepackage{amsmath}
\usepackage{amsthm}
\usepackage{amssymb}

\usepackage{float}
\usepackage{hyperref}
\usepackage{fullpage}
\usepackage{graphicx}
\usepackage{caption}
\usepackage{subcaption}

\usepackage{color}
\usepackage{hyperref}

\newtheorem{theo}{Theorem}

\newtheorem{deff}{Definition}

\title{A Note on a Recent Conjecture in Mathematics Magazine}

\author[Malinovsky]{Yaakov Malinovsky}
\address[Y.~Malinovsky]{Department of Mathematics and Statistics,
University of Maryland, Baltimore County, Baltimore, MD  21250}
\email{\href{mailto:yaakovm@umbc.edu}{yaakovm@umbc.edu}}

\begin{document}
\maketitle
David Foster and T\"{urkay} Yolcu, in the final section of their paper (Math. Mag. 98 (2025), no. 2)\cite{FY2025}, proposed an open problem, which can be stated as follows.
Suppose we perform an experiment by randomly and independently dropping, one by one, $N$ balls into $n$ boxes labeled $1, 2, \ldots, n$, such that each ball lands in box $i$ with probability $p_i$, for $i = 1, \ldots, n$, where $p_i > 0$ and $\sum_{i=1}^{n} p_i = 1$. Let $X$ be the number of boxes that contain at least one ball after the experiment is completed, and let $E_p$ denote the expectation of $X$, where $p = (p_1, p_2, \ldots, p_n)$. Without loss of generality, we may assume that $p_1 \geq p_2 \geq \cdots \geq p_n$. David Foster and T\"{urkay} Yolcu conjectured that
\begin{equation}
\label{ex:Conj}
\underset{p \in \left\{p:\,\,\, p_i > 0,\,\, \sum_{i=1}^{n} p_i = 1\right\}}{\operatorname{argmax}}  E_p=\left(1/n,\ldots,1/n\right),
\end{equation}
i.e., the expectation is maximal under the uniform distribution.
In fact, this conjecture is true and was proven in \cite{WY1973} using the famous theorem of Issai Schur \cite{S1923}. The technique is easy to explain and apply, and it has been used in many problems of a similar type. For example, see \cite{CW1991} and \cite{MOA2009}. In the following few lines, we outline the proof from \cite{WY1973}, using formal definitions from \cite{MOP1967}, and connect this problem to other related results in the literature.
\begin{deff}[\cite{MOP1967}]
If $a_1\geq \cdots \geq a_n$,\,\,$b_1\geq \cdots \geq b_n$, $\sum_{j=1}^{k}a_j\geq \sum_{j=1}^{k}b_j$ for $k=1, 2, \ldots, n-1$, and $\sum_{j=1}^{n}a_j=\sum_{j=1}^{n}b_j$,
then $a=\left(a_1,\ldots,a_n\right)$ is said to majorize $b=\left(b_1,\ldots,b_n\right)$, written $a\succ b$.
\end{deff}
\begin{deff}[\cite{MOP1967}]
A real function $\varphi$ of $n$ real variables $(x_1,\ldots x_n)=x$ is said to be a $\it Schur$ function if for very pair $i\neq j$, $ (x_i-x_j)\left(\frac{\partial{\varphi(x)}}{\partial{x_i}}-\frac{\partial{\varphi(x)}}{\partial{x_j}}\right)\geq 0.$
\end{deff}

Also, \cite{MOP1967} presents the following theorem, originally proved by Schur for the case where $x_i\in (0,\infty)$, and later extended by Ostrowski to any open subset of the real line (see also \cite{MOA2009}).
\begin{theo}[Schur \cite{S1923}, Ostrowski \cite{O1952}]
\label{Schur}
Let $\varphi(x)$ be defined for $x_1\geq \cdots \geq x_n$. Then $\varphi(a)\geq \varphi(b)$  for all $a\succ b$ if and only if $\varphi$ is a {\it Schur} function.
\end{theo}

Wong and Yue \cite{WY1973} showed that $E_p=n-\sum_{j=1}^{n}(1-p_i)^N.$
Since $\varphi(p)=\sum_{j=1}^{n}(1-p_i)^N$ is a Schur function, by appealing to Theorem~\ref{Schur}, one obtains:
\begin{equation}
\label{SO}
\varphi(q)\leq \varphi(p)\,\,\,\, \text {if and only if}\,\,\,\, (q_1,\ldots,q_n)=q\prec p=(p_1,\ldots,p_n).
\end{equation}

In particular, \eqref{SO} holds for $(q_1,\ldots,q_n)=(1/n,\ldots,1/n)$, since for any $p = (p_1, p_2, \ldots, p_n)$, $(1/n,\ldots,1/n)\prec p$.
Therefore, the conjecture \eqref{ex:Conj} holds.

In fact, Wong and Yue \cite{WY1973} proved a much more general result. Namely, they showed that for any $q\prec p$ and for all $k$,
$$
P_{q}\left(X\leq k\right)\leq P_{p}\left(X\leq k\right).
$$
This result also implies the conjecture \eqref{ex:Conj}, since ${\displaystyle E_{p}=\sum_{k=1}^{n}P_{p}\left(X\geq k\right)}$.

Finally, it is worth mentioning that the problem of empty boxes, defined as $n-X$ in the described experiment, has been extensively studied. See, for example, the book by Kolchin, Sevast'yanov, and Chistyakov \cite{KSC1978}, the references therein, and the papers that cite this book.

\end{document}